\documentclass[11 pt,leqno]{article} 

\usepackage[all]{xy}
\usepackage{amsfonts}
\usepackage{mathrsfs}
\usepackage{amsmath}
\usepackage{amssymb}
\usepackage{amsthm}
\usepackage{latexsym}
\usepackage{epsfig}
\usepackage{verbatim}
\title{Extending Johnson's and Morita's homomorphisms to the mapping class group}

\author{Matthew B Day}

\date{24 September 2007
\\Originally published in: Algebraic \& Geometric Topology 7 (2007) 1297--1326
\\{\tt http://msp.warwick.ac.uk/agt/2007/07/p050.xhtml}}

\theoremstyle{plain}
\newtheorem{theorem}{Theorem}[section]
\newtheorem{proposition}[theorem]{Proposition}
\newtheorem{lemma}[theorem]{Lemma}
\newtheorem{corollary}[theorem]{Corollary}
\newtheorem{claim}[theorem]{Claim}
\newtheorem{maintheorem}{Main Theorem}
\newtheorem{maincorollary}[maintheorem]{Main Corollary}

\theoremstyle{definition} 
\newtheorem*{remark}{Remark}
\newtheorem{definition}[theorem]{Definition}

\numberwithin{equation}{section}

\newcommand\R{\mbox{$\mathbb{R}$}}
\newcommand\C{\mbox{$\mathbb{C}$}}
\newcommand\Z{\mbox{$\mathbb{Z}$}}
\newcommand\Q{\mbox{$\mathbb{Q}$}}

\def\Hom{\mathop{\mathrm{Hom}}\nolimits}
\def\Diff{\mathop{\mathrm{Diff}}\nolimits}
\def\Aut{\mathop{\mathrm{Aut}}\nolimits}
\def\Sp{\mathop{\mathrm{Sp}}\nolimits}
\newcommand\into\hookrightarrow
\newcommand\isomarrow{\hspace{0.7 em}{}^{\cong}\hspace{-1.3 em}\longrightarrow}
\def\co{\colon\thinspace}
\newcommand\bdI[1]{\mathcal{I}(#1)}
\newcommand\bdIg[1]{\mathit{I}_{g,1}(#1)}
\newcommand\bdM{\mathit{Mod}}
\newcommand\bdMg{\mathit{Mod}_{g,1}}
\newcommand\mohom[1]{\tilde\tau_#1}
\newcommand\mohomg[1]{\tilde\tau_{g,1}(#1)}
\newcommand\johom[1]{\tau_#1}
\newcommand\johomg[1]{\tau_{g,1}(#1)}
\newcommand\lv[1]{\mathcal{L}_{#1}}
\newcommand\lvr[1]{\lv{#1}\otimes\R}
\newcommand\gam[1]{\Gamma_{#1}}
\newcommand\lag[1]{\mathfrak{g}_{#1}}
\newcommand\lal[1]{\mathfrak{l}_{#1}}
\newcommand\tomohom[1]{\overline{\epsilon}_{#1}}
\newcommand\emomap[1]{\epsilon_{#1}}
\newcommand\CmodB[1]{C_3(\lag{#1})/B_3(\lag{#1})}
\newcommand\exdiff{\widetilde{d^2}}
\begin{document}

\maketitle

\begin{abstract}
We extend certain homomorphisms defined on the higher Torelli subgroups of the mapping class group to crossed homomorphisms defined on the entire mapping class group.
In particular, for every $k\geq 2$, we construct a crossed homomorphism $\epsilon_k$ which extends Morita's homomorphism $\tilde \tau_k$ to the entire mapping class group.
From this crossed homomorphism we also obtain a crossed homomorphism extending the $k$th Johnson homomorphism $\tau_k$ to the mapping class group.

D Johnson and S Morita obtained their respective homomorphisms by considering the action of the mapping class group on the nilpotent truncations of the surface group; our approach is to mimic Morita's construction topologically by using nilmanifolds associated to these truncations.
This allows us to take the ranges of these crossed homomorphisms to be certain finite-dimensional real vector spaces associated to these nilmanifolds.
\end{abstract}

%%%%%%%%%%%%%%%%%%%%   Start of main body of article

\section{Introduction}\label{se:introduction}
Let $\Sigma=\Sigma_{g,1}$ be a compact, connected, oriented smooth surface of genus $g\geq 3$ with a single boundary component and a basepoint $*$ on the boundary component.
The mapping class group $\bdMg$ of $\Sigma$ is the group $\Diff(\Sigma,\partial\Sigma)$ of diffeomorphisms of $\Sigma$ fixing $\partial\Sigma$ pointwise, modulo isotopies fixing $\partial\Sigma$ pointwise.

Let $\pi=\pi_1(\Sigma,*)$ be the fundamental group of $\Sigma$, which is a free group on $2g$ generators.
The obvious action of $\Diff(\Sigma,\partial\Sigma)$ on the based loops of $\Sigma$ induces an action of $\bdMg$ on $\pi$.
Consider the lower central series of $\pi$, given by $\pi^{(0)}=\pi$ and $\pi^{(k+1)}=[\pi,\pi^{(k)}]$.
The \emph{$k$th nilpotent truncation of $\pi$} is the group $\gam{k}=\pi/\pi^{(k-1)}$.
This is a nilpotent group of class $(k-1)$.
Since the group $\pi^{(k)}$ is characteristic in $\pi$, the action of $\bdMg$ on $\pi$ descends to an action on $\gam{k}$ for each $k$.
So we have a representation $\bdMg\to\Aut(\gam{k})$.
The \emph{$k$th Torelli group} $\bdIg{k}$ is the kernel of this representation.
These groups have been studied by Johnson and others (see~\cite{Johnson83}).
Andreadakis defined analogues of these groups in $\Aut(\pi)$ in~\cite{Andreadakis}. 

Let $H=H_1(\Sigma;\Z)$ be the integral first homology group of the surface and let $\lv{k+1}$ be the abelian group $\pi^{(k-1)}/\pi^{(k)}$.
There is an exact sequence (see Morita~\cite[Proposition~2.3]{Morita93Invent}):
\[0\to \Hom(H,\lv{k+1})\to \Aut(\gam{k+1})\to \Aut(\gam{k})\to 1\]
The image of $\bdIg{k}$ under the representation $\bdMg\to\Aut(\gam{k+1})$ lies in the image of $\Hom(H,\lv{k+1})\hookrightarrow \Aut(\gam{k+1})$.
So by restriction, we have a homomorphism $\johomg{k}$ for any $k\geq 2$:
\[\johomg{k}\co\bdIg{k}\to\Hom(H,\lv{k+1})\]
This is called the \emph{$k$th Johnson homomorphism}.
It was defined by Johnson in~\cite{Johnson80}.

Consider the standard (bar) group cohomology chains $C_*(\pi)$ of $\pi$.
The action of $\bdMg$ on $\pi$ induces an action on $C_*(\pi)$.
The abelian group $C_1(\pi)$ is the free abelian group with generating set $\pi$, so we have an element $[\ell]\in C_1(\pi)$, where $\ell\in \pi$ is the class of the boundary loop in $\Sigma$.
We denote the normal closure of $\langle\ell\rangle$ in $\pi$ by $\overline{\langle\ell\rangle}$.
Select a $2$--chain $C\in C_2(\pi)$ such that $\partial C=-[\ell]$ and the image of $C$ in $C_2(\pi/\overline{\langle\ell\rangle})$ is a cycle 
representing the fundamental class of $\pi/\overline{\langle\ell\rangle}$.
For $[\phi]\in\bdMg$, the chain 
$[\phi]\cdot C- C$
is in $Z_2(\pi)$ because $[\ell]$ is fixed by $[\phi]$, and it is in $B_2(\pi)$ because $H_2(\pi)=0$.
Pick a chain $D\in C_3(\pi)$ that bounds $[\phi]\cdot C-C$.
Suppose further that $[\phi]\in\bdIg{k}$.
Then $[\phi]$ fixes the image of $C$ in $C_3(\gam{k})$ and the image of $D$ in $C_3(\gam{k})$ is a cycle.
The \emph{$k$th Morita homomorphism} $\mohomg{k}$ is the map taking $[\phi]$ to the homology class of this cycle:
\[\mohomg{k}\co\bdIg{k}\to H_3(\gam{k})\]
Morita proved that this is a well-defined homomorphism in~\cite[Theorem~3.1]{Morita93}, where he also proved that it refines $\johomg{k}$ (see Theorem~\ref{th:MoritaTheorem} below).
When the genus is suppressed, we use the notation $\johom{k}$ for $\johomg{k}$  and $\mohom{k}$ for $\mohomg{k}$. 

Since $\pi$ is a free group, each $\gam{k}$ is a finitely generated, torsion-free nilpotent group, and is therefore, by a theorem of Mal'cev, a lattice in a unique simply-connected nilpotent Lie group $G_k$ over the reals.
Let $X_k=G_k/\gam{k}$, the homogeneous space of $\gam{k}$, and let $\lag{k}$ be the Lie algebra of $G_k$.
The standard Lie algebra chain complex $C_*(\lag{k})$ on $\lag{k}$ is a chain complex with chains $C_n(\lag{k})$ equal to the exterior power $\Lambda^n\lag{k}$ and a boundary determined by the Lie bracket.
Note that the vector spaces $C_n(\lag{k})$ are finite-dimensional for every $n$ and $k$.
By a theorem of Nomizu~\cite{Nomizu}, the homology $H_*(\lag{k})$ of this complex is isomorphic to $H_*(X_k;\R)$.
This is explained in Section~\ref{se:nilmanifolds} below.

A \emph{crossed homomorphism} from a group $G$ to a $G$--module $M$ is a map $G\to M$ satisfying a product rule that is twisted by the action of $G$ on $M$:
\[f(gh)=g\cdot f(h)+f(g)\]
This is the same thing as a group cohomology $1$--cocycle in $Z^1(G;M)$, so we may speak of the cohomology class of a crossed homomorphism in $H^1(G;M)$.

The issue of extending Johnson's homomorphisms is first addressed in Morita's paper~\cite{Morita93Invent}, where he defines a crossed homomorphism extending $\johomg{2}$.
Later, Morita constructed an extension of $\johomg{3}$ in~\cite{Morita96}, and Perron in~\cite{Perron} built an extension of $\johomg{3}$ by different methods.
Hain showed in~\cite[Section~14.6]{Hain97} that $\johomg{k}$ extends to a representation of $\bdMg$ in a semi-direct product of a unipotent algebraic group over $\Q$ with $\Sp_{2g}(\Q)$.
We note that while semi-direct products and crossed homomorphisms are related, the results of this paper do not appear to follow directly from Hain's results.
Kawazumi constructed a different class of maps in~\cite{Kawazumi05} that extend the higher Johnson homomorphisms to $\Aut(\pi)$, and showed in~\cite{Kawazumi06} that these maps are not crossed homomorphisms.

\begin{maintheorem}\label{mt:ExtMoritaHom}
Fix $k\geq 2$.
There is an action of $\bdMg$ on the finite-dimensional real vector space $\CmodB{k}$, and with respect to this action, there is a crossed homomorphism 
\[\emomap{k}\co\bdMg\to \CmodB{k}\]
that extends the Morita homomorphism:
\[
\xymatrix{
\bdIg{k}\ar@{^{(}->}[d]\ar[rr]^{\mohomg{k}} && H_3(\gam{k})\ar@{^{(}->}[d]\\
\bdMg\ar[rr]^{\emomap{k}} && \CmodB{k}
}
\]
\end{maintheorem}
The unlabeled map on the right is the composition: 
\[H_3(\gam{k})\isomarrow H_3(X_k;\Z)\into H_3(X_k;\R)\isomarrow H_3(\lag{k}) \into \CmodB{k}\]

We construct this crossed homomorphism in Section~\ref{ss:CoreConstruction} below.
The construction is a generalization of a topological construction for $\johomg{2}$ (see Johnson~\cite{Johnson80} or Hain~\cite{Hain95}) and is similar to a different topological construction of $\mohomg{k}$ due to A. Heap (see \cite[Theorem~4]{Heap} and \cite[Theorem~22]{Heap}).

The projection $\gam{k+1}\to\gam{k}$ induces $G_{k+1}\to G_k$ on Mal'cev completions, which then induces $\lag{k+1}\to\lag{k}$ on Lie algebras.
Let $\lal{k+1}$ be the kernel of $\lag{k+1}\to\lag{k}$, and let $\lag{k+1}^{(1)}$ denote the commutator subalgebra of $\lag{k+1}$.
From Main Theorem~\ref{mt:ExtMoritaHom}, we obtain
\begin{maintheorem}\label{mt:ExtJohnsonHom}
Fix $k\geq 2$.
There is a $\bdMg$--equivariant map 
\[\exdiff\co\CmodB{k}\to C_2(\lag{k+1})/(\lag{k+1}^{(1)}\wedge\lal{k+1})\]
 such that the crossed homomorphism
\[\exdiff\circ \emomap{k}\co\bdMg\to C_2(\lag{k+1})/(\lag{k+1}^{(1)}\wedge\lal{k+1})\]
extends the Johnson homomorphism:
\[
\xymatrix{
\bdIg{k}\ar@{^{(}->}[d]\ar[rrr]^{\johomg{k}} &&& \Hom(H,\lv{k+1})\ar@{^{(}->}[d]\\
\bdMg\ar[rrr]^{\exdiff\circ\emomap{k}} &&& C_2(\lag{k+1})/(\lag{k+1}^{(1)}\wedge\lal{k+1}) \\
}
\]
\end{maintheorem}
The inclusion on the right is explained in the proof of the theorem.

\begin{remark}
It is a fact of group cohomology theory (see Brown~\cite[chapter III.5]{Brown}) that a crossed homomorphism on a subgroup may be extended to the entire group by expanding the target module to a larger module, called a co-induced module.
As an abelian group, the co-induced module of $H_3(\gam{k})$ from $\bdIg{k}$ to $\bdMg$ is 
a direct sum of copies of $H_3(\gam{k})$ indexed over the cosets of $\bdIg{k}$ in $\bdMg$.
Our approach finds a crossed homomorphism extending $\mohomg{k}$ with range in a finite-dimensional vector space, while the co-induced module approach would make the range an infinite-rank \mbox{$\bdMg$-module} with no additional structure.
\end{remark}

\begin{remark}
As is common for crossed homomorphisms, the definitions of these crossed homomorphisms involve an essential choice.
Varying this choice changes their values on mapping classes outside of $\bdIg{k}$, but it does not change the cohomology classes 
\begin{gather*}
[\emomap{k}]\in H^1(\bdMg;\CmodB{k})\\
\tag*{\text{and}}
[\exdiff\circ\emomap{k}]\in H^1(\bdMg;C_2(\lag{k+1})/(\lag{k+1}^{(1)}\wedge\lal{k+1})).
\end{gather*}
This is explained in Proposition~\ref{pr:lostbasemap} below.
\end{remark}

As an application of Main Theorem~\ref{mt:ExtMoritaHom} and Main Theorem~\ref{mt:ExtJohnsonHom}, we obtain the following Corollary:
\begin{maincorollary}\label{mc:sdproduct}
We have the following injections for all $k\geq 2$:
\begin{align*}
\bdMg/\bdIg{2k-1}&\into(\bdMg/\bdIg{k})\ltimes (\CmodB{k})\\
\tag*{\text{by}} \alpha\bdIg{2k-1}&\mapsto(\alpha\bdIg{k},\emomap{k}(\alpha))\\
\tag*{\text{and}}
\bdMg/\bdIg{k+1}&\into(\bdMg/\bdIg{k})\ltimes (C_2(\lag{k+1})/(\lag{k+1}^{(1)}\wedge\lal{k+1}))\\
\tag*{\text{by}}\alpha\bdIg{k+1}&\mapsto(\alpha\bdIg{k},\exdiff\circ\emomap{k}(\alpha))
\end{align*}
The semi-direct products are taken with respect to actions of $\bdMg/\bdIg{k}$ that are induced from the usual actions of $\bdMg$ on the respective vector spaces.
\end{maincorollary}

Section~\ref{se:nilmanifolds} develops the ideas from the theory of Nilpotent Lie groups that are relevant to our construction.  
Section~\ref{se:mainconstruction} contains the construction of $\emomap{k}$ and Section~\ref{se:HomomorphismEquivalence} completes the proof of Main Theorem~\ref{mt:ExtMoritaHom}.
Main Theorem~\ref{mt:ExtJohnsonHom} and Main Corollary~\ref{mc:sdproduct} are proven in Section~\ref{se:proofof}.

\subsection{Acknowledgements}\label{ss:acknowledgements}
This research was supported by a graduate research fellowship from the National Science Foundation.
My thanks go out to the following people for their help with this project:
my thesis advisor Benson Farb for his constant guidance, 
David Witte Morris for conversations that led to some of the ideas in this paper,
Robert Young for conversations and for indicating a simplification in the proof that $\emomap{k}$ is well-defined,
and Richard Hain, Nariya Kawazumi, Nathan Broaddus, Asaf Hadari and T. Andrew Putman for comments on earlier versions of this paper.
I would also like to thank the anonymous referee for many useful comments and suggestions.

%%%%%%%%%%%%%%%%%%%%%%%%%%%%%%%%%%%%%%%%%%%%%%%%%%%%%%%%%%%%%%%%%%%%%%%%%%%%%%
\section{Preliminaries about nilmanifolds}\label{se:nilmanifolds}
This section develops some ideas about nilmanifolds, including
 some consequences of Mal'cev's rigidity theory and some notions relating to the homology of compact homogeneous spaces of nilpotent Lie groups.

The $3$--dimensional Heisenberg group is a good example to keep in mind throughout this section.

%%%%%%%%%%%%%%%%%%%%%%%%%%%%%%%%%%%%%%%%%%%%%%%%
\subsection{Homology of Lie algebras and nilpotent homogeneous spaces}
\label{ss:Liealgebrahomology}
One of the key ingredients in the construction of $\emomap{k}$ is the homology of Lie algebras.
Koszul~\cite{Koszul} is a classical reference for this theory; Weibel~\cite[chapter 7]{Weibel} is a modern one.
All vector spaces and Lie algebras in this paper are defined over the reals.

Suppose that $G$ is a Lie group with Lie algebra $\lag{}$.  
Let $X=G/\Gamma$ be the homogeneous space of $G$ with respect to a lattice $\Gamma$ in $G$.
Let $\lag{}^*$ be the dual vector space to $\lag{}$, \emph{ie} the vector space of linear functionals on $\lag{}$.
We use these to define the following graded vector spaces:
\[C_*(\lag{})=\bigoplus_n\Lambda^n\lag{},\quad\mbox{and}\quad C^*(\lag{})=\bigoplus_n\Lambda^n\lag{}^*\]
The vector space $C_*(\lag{})$ has a degree-one boundary map $\partial$ given by:
\begin{equation}
\begin{split}
\partial (V_1\wedge&\ldots \wedge V_n)=\\
&\sum_{1\le i<j\le n}(-1)^{i+j+1}[V_i,V_j]\wedge V_1\wedge\ldots\wedge \widehat{V}_i\wedge\ldots\wedge\widehat{V}_j\wedge\ldots\wedge V_n
\end{split}
\label{eq:KoszulBoundary}
\end{equation}
where $V_1,\ldots, V_n\in \lag{}$ and where the caret $(\widehat{\phantom{V}})$ indicates omission.

By definition, we have an evaluation pairing between $C_n(\lag{})$ and $C^n(\lag{})$ for any $n$.
The degree $-1$ differential $d$ on $C^*(\lag{})$ is defined by
\[d\alpha(\overline{V})=\alpha(\partial \overline{V})\]
where $\alpha\in C^n(\lag{})$ and $\overline{V}\in C_n(\lag{})$ for some $n$.

It is a straightforward consequence of the Jacobi identity that
the pair $(C_*(\lag{}),\partial)$ is a chain complex and the pair $(C^*(\lag{}),d)$ is a cochain complex.
A proof can be found in Koszul~\cite{Koszul}.

We use left-propagation to connect Lie algebras to Lie groups and homogeneous spaces. 
Let $C^*(G;\R)$ and $C^*(X;\R)$ be the De Rham cochain complexes on $G$ and $X$ respectively.

\begin{definition}
For $\alpha\in C^*(\lag{})$, define $L(\alpha)$ to be the unique differential form on $G$ that is invariant under the left action of $G$ and restricts to $\alpha$ at the identity $e$. 
This defines a map $L:C^*(\lag{})\to C^*(G;\R)$ called the \emph{left-propagation map}.

For $\alpha\in C^*(\lag{})$, the differential form $L(\alpha)$ descends to a form on $X$, so we also get a left propagation map $L:C^*(\lag{})\to C^*(X;\R)$. 
\end{definition}

One can write down an explicit formula for the value of $L(\alpha)$ at a point $g$ in $G$; from this we deduce that $L$ is well-defined.

The following theorem of Nomizu is from~\cite{Nomizu}.
\begin{theorem}[Nomizu]\label{th:Nomizu}
If $X=G/\Gamma$ is the homogeneous space of a connected, simply-connected nilpotent Lie group $G$ with respect to a lattice $\Gamma$ and $\lag{}$ is the Lie algebra of $G$, then the left-propagation map
\[L\co C^*(\lag{})\into C^*(X;\R)\]
induces an isomorphism
\[L_*\co H^*(\lag{})\isomarrow H^*(X;\R)\]
on cohomology groups.
\end{theorem}

The goal of this subsection is to construct an adjoint to the map $L_*$.
Let $C_*(X)$ and $C_*(X;\R)$ denote the piecewise-smooth singular homology chains of $X$ with coefficients in $\Z$ and $\R$ respectively.

\begin{definition}\label{de:vectorchainmap}
For each $n$, the vector space $C_n(\lag{})$ is finite dimensional and is canonically identified with its double dual, $(C^n(\lag{}))^*$.
We can therefore think of the transpose to $L$ as a map:
\[L^*\co (C^*(X;\R))^*\to C_*(\lag{})\] 

Define 
\[v\co C_*(X;\R)\to C_*(\lag{})\] 
to be the composition of $L^*$ with the map
$C_*(X;\R)\to (C^*(X;\R))^*$
given by:
\[C\mapsto (\alpha\mapsto \int_C\alpha)\]
\end{definition}
We can think of $v$ as taking singular chains on $X$ and smoothing them out into left-invariant multi-vector fields.

\begin{claim}
The map $v$ is a chain complex map.
\end{claim}

\begin{proof}
We only need to prove that $v$ respects the boundary maps.
The map $L^*$ obviously respects the boundary, and the integration map respects the boundary by Stokes's theorem. 
\end{proof}

\begin{lemma} \label{le:ObviousDuality}
For each $n$, for every $\alpha\in C^n(\lag{})$ and $C\in C_n(X;\R)$, we have:
\begin{equation}\alpha(v(C))=\int_C L(\alpha)\label{eq:ObviousDuality}\end{equation}
\end{lemma}

\begin{proof}This is immediate from the definition.\end{proof}

\begin{lemma}\label{le:CoordinateForm}
Suppose we have a basis $B$ for $C_n(\lag{})$.
This basis determines a dual basis $B^*$ for $C^n(\lag{})$ and a map $()^*:B\to B^*$.
This gives us the following coordinate formula for $v_n:$
\[v_n(C)=\sum_{\overline{V}\in B}\big(\int_C \overline{V}^*\big)\overline{V}\]
where $C\in C_n(X;\R)$.
\end{lemma}
\begin{proof}
This is also immediate from the definition.
\end{proof}

\begin{corollary}\label{co:dualofNomizuTheorem}
Suppose $X=G/\Gamma$ is a compact homogeneous space of a connected, simply-connected nilpotent Lie group $G$ with respect to a discrete subgroup $\Gamma$ and $\lag{}$ is the Lie algebra of $\Gamma$.
Then the map \[v\co C_*(X;\R)\to C_*(\lag{})\] induces an isomorphism 
\[v_*\co H_*(X;\R)\isomarrow H_*(\lag{})\]
on homology groups.
\end{corollary}
\begin{proof}
Fix an $n$.
For $\alpha\in Z^n(\lag{})$ and $C\in Z_n(X;\R)$, Equation~(\ref{eq:ObviousDuality}) applies and induces the following evaluation identity on homology and cohomology:
\[[\alpha]([v(C)])=[L(\alpha)]([C])\]

This means that with respect to the evaluation pairings, $v_*$ is adjoint to $L_*$.
The hypotheses are the same as those for Theorem~\ref{th:Nomizu}, so we have that $L_*$ is an isomorphism.
The adjoint of an isomorphism is an isomorphism, so $v_*$ is also an isomorphism.
\end{proof}
%%%%%%%%%%%%%%%%%%%%%%%%%%%%%%%%%%%%%%%%%%%%%%%%%%%%%%%%%%%%%%%%%%%%%%%%%%%5555
\subsection{Actions on nilmanifolds}
\label{ss:actionsonX}
The starting point of this subsection is Mal'cev rigidity theory, for which Raghunathan~\cite[chapter II]{Raghunathan} is a good reference. 
Throughout this subsection, $\Gamma$ is finitely generated, torsion-free nilpotent group.
The existence statement of Theorem~\ref{th:Mal'cevCompletion} follows from~\cite[Theorem~2.18]{Raghunathan} and the uniqueness statement of Theorem~\ref{th:Mal'cevCompletion} and Theorem~\ref{th:CompleteMaps} follow from \cite[Theorem~2.11]{Raghunathan}.

\begin{theorem}[Mal'cev]\label{th:Mal'cevCompletion}
There is a unique connected, simply-connected, nilpotent Lie group $G$ such that $\Gamma\into G$ as a lattice.
This Lie group is called the \emph{Mal'cev completion} of $\Gamma$.
\end{theorem}

\begin{theorem}[Mal'cev]\label{th:CompleteMaps}
Every map $\Gamma\to\Gamma'$ of torsion-free, nilpotent groups extends to a unique map of $G\to G'$ of their Mal'cev completions.
In particular, if $\Gamma'=\Gamma$, this describes an injective homomorphism $\Aut(\Gamma)\into \Aut(G)$.
\end{theorem}

Let $G$ be the Mal'cev completion of $\Gamma$ and consider the homogeneous space $X=G/\Gamma$.
This is a manifold with the coset of the identity $[e]$ as a natural basepoint.
Since torsion-free nilpotent Lie groups are aspherical, the space $X$ is a $K(\Gamma,1)$.
\begin{proposition}\label{pr:actiononX}
There is an action of $\Aut(\Gamma)$ on $(X,[e])$, such that the induced action on the fundamental group is the usual action of $\Aut(\Gamma)$ on $\Gamma$.
\end{proposition}

\begin{proof}
This follows immediately from Theorem~\ref{th:CompleteMaps}.
\end{proof}

The action of $\Aut(\Gamma)$ on $G$ induces an action on $\lag{}$ by taking derivatives at the identity. 
The actions on $\lag{}$ and $C_*(\lag{})$ are denoted with push-forward notation (\emph{e.g.} $\phi_*\overline{V}$).
The right actions on their duals are denoted with pull-back notation on the left (\emph{e.g.} $\phi^*\alpha$).

\begin{claim}
The left-propagation map $L\co C^*(\lag{})\to C^*(X;\R)$ is $\Aut(\Gamma)$--equivariant.
\end{claim}

\begin{proof}
One can show by a computation that the map
\begin{align*}
C^*(X;\R)&\to \lag{}^*\\
\tag*{\text{by}}
\alpha&\mapsto \alpha|_{[e]}
\end{align*}
is equivariant.
The claim then follows from the definition of $L$.
\end{proof}

\begin{claim}
The map $v\co C_*(X;\R)\to C_*(\lag{})$ is $\Aut(\Gamma)$--equivariant.
\end{claim}

\begin{proof}
This is immediate from the fact that $L$ is equivariant and Lemma~\ref{le:ObviousDuality}.
\end{proof}

\begin{proposition} \label{pr:equivariantsection}
There is an $\Aut(\Gamma)$--equivariant map 
\[s\co \Gamma\to C^\infty((S^1,*),(X,[e]))\]
that is a section to the projection $C^\infty((S^1,*),(X,[e]))\to \pi_1(X)=\Gamma$.
\end{proposition}
Here $C^\infty((S^1,*),(X,[e]))$ is the set of smooth, based loops of $X$.
\begin{proof}
For $\gamma\in\Gamma$,
let $f(\gamma)\co \R\to G$ be the one-parameter family in $G$ such that $f(\gamma)(1)=\gamma$.
This exists and is unique by Theorem~\ref{th:CompleteMaps}.
The map $f(\gamma)$ descends to a map $s(\gamma)\co (S^1,*)\to (X,[e])$.
The map $s(\gamma)$ is obviously smooth, giving us 
$s\co \Gamma\to C^\infty((S^1,*),(X,[e]))$.
Let $\rho\co \Gamma\to\Aut(G)$ denote the action of $\Aut(\Gamma)$ on $G$. 
If $\phi\in\Aut(\Gamma)$, then $\rho(\phi)\circ f(\gamma)$ is a one-parameter family with $\rho(\phi)\circ f(\gamma)(1)=\phi\cdot\gamma$.
By uniqueness, $\rho(\phi)\circ f(\gamma)=f(\phi\cdot\gamma)$, so $s$ is $\Aut(\Gamma)$--equivariant.
The fact that $s$ is a section to  $C^\infty((S^1,*),(X,[e]))\to\Gamma$ follows immediately from Proposition~\ref{pr:actiononX}.
\end{proof}
%%%%%%%%%%%%%%%%%%%%%%%%%%%%%%%%%%%%%%%%%%%%%%%%%%%%%%%%%%%%%%%%%%%%%%%%555
\section{The main construction}\label{se:mainconstruction}
\subsection{Preliminary definitions}\label{ss:maindefinitions}
Fix a $k\geq 2$.  Since $g$ is forever fixed and no comparisons to surfaces of other genera will be made, we will drop it from the notation.
We recall the definitions of the following objects, which were mentioned in the introduction:
the group $\gam{k}=\pi/\pi^{(k-1)}$ is the class $(k-1)$ nilpotent truncation of the surface group $\pi$,
the Lie group $G_k$ is the Mal'cev completion of $\gam{k}$,
the Lie algebra $\lag{k}$ is the tangent space at the identity of $G_k$, and
the space $X_k$ is the homogeneous space $G_k/\gam{k}$.
Also recall that we have an action of $\bdM$ on $\gam{k}$.

We obtain actions of $\bdM$ on $G_k$, $\lag{k}$ and $X_k$ via $\bdM\to\Aut(\gam{k})$, since $\Aut(\gam{k})$ acts on these same objects by Theorem~\ref{th:CompleteMaps} and Proposition~\ref{pr:actiononX}.
The induced action of $\bdM$ on $C_3(\lag{k})/B_3(\lag{k})$ gives it the module structure we use to prove Main Theorem~\ref{mt:ExtMoritaHom}. 

To continue, we need to define a space of maps that relate $\Sigma$ and $X_k$.
As the first step, we need to find an analogue of the boundary loop $\partial\Sigma$ in the space $X_k$.
\begin{proposition}\label{pr:loopinX}
If $p\co \pi\to\gam{k}$ is the canonical projection, and $[\partial\Sigma]\in\pi_1(\Sigma,*)$ is the boundary loop class, then
there is a smooth representative $\ell_k$ of $p_*[\partial\Sigma]\in\pi_1(X_k,[e])$ that is pointwise fixed by the action of $\bdM$ on $X_k$.
\end{proposition}

\begin{proof}
By Proposition~\ref{pr:equivariantsection}, we have a section $s$ from $\gam{k}$ to the space of based smooth loops in $X_k$, such that the image of $s$ is invariant under the action of $\bdM$.
Since the boundary loop class $[\partial\Sigma]$ is fixed by the action of $\bdM$, so is its image $p_*[\partial\Sigma]$ in $\gam{k}$.
So if we let $\ell_k$ be $s(p_*[\partial\Sigma])$, then $\ell_k$ is pointwise fixed.
\end{proof}

Let $\ell_k$ be as defined by Proposition~\ref{pr:loopinX}.
We also use the symbol $\ell_k$ to refer to the image $\ell_k(S^1)$.
Fix a piecewise-smooth map \[i_{\partial\Sigma}\co (\partial\Sigma,*)\to(\ell_k,[e])\] that sends $[\partial\Sigma]$ to $[\ell_k]$ at the fundamental group level.

\begin{definition}\label{de:spaceofmaps}
Define $S_k$ to be the space of all piecewise-smooth maps $\Sigma\to X_k$ that restrict to $i_{\partial\Sigma}$ on $\partial\Sigma$ and induce the canonical projection $\pi\to\gam{k}$ on fundamental groups.
We topologize $S_k$ by giving it the compact--open topology.

For $\phi\in\Diff(\Sigma,\partial\Sigma)$ and $i\in S_k$, we define $\phi\cdot i$ by the formula 
\[(\phi\cdot i)(p)=\phi\cdot (i(\phi^{-1}(p)))\]
for $p\in\Sigma$.
Here $\Diff(\Sigma,\partial \Sigma)$ acts on $X_k$ through its projection to $\bdM$.

If $F\co \Sigma\times [0,1]\to X_k$ is a homotopy through maps in $S_k$, we define $\phi\cdot F$ by
\[(\phi\cdot F)_t=\phi\cdot F_t\]
for $t\in[0,1]$.
\end{definition}

\begin{claim}
The action of $\Diff(\Sigma,\partial\Sigma)$ on $S_k$ is well-defined.
\end{claim}

\begin{proof}
We only need to show that for $i\in S_k$ and $\phi\in\Diff(\Sigma,\partial\Sigma)$, we have $\phi\cdot i\in S_k$.
It follows from Proposition~\ref{pr:actiononX} that $\phi\cdot i$ induces the canonical projection on fundamental groups, and it follows from Proposition~\ref{pr:loopinX} that $(\phi\cdot i)|_{\partial \Sigma}=i_{\partial\Sigma}$.
\end{proof}

\begin{lemma}\label{le:Snonempty}
The space $S_k$ is nonempty.
\end{lemma}

\begin{proof}
Let $h\co \Sigma\to X_k$ be a smooth map inducing the canonical projection on fundamental groups.
Such a map exists because $\Sigma$ and $X_k$ are both smooth, aspherical manifolds.

Since $h|_{\partial X}$ and $i_{\partial\Sigma}$ induce the same map on fundamental groups, we have a homotopy $F\co \partial \Sigma\times [0,1]\to X_k$ from $i_{\partial\Sigma}$ to $h|_{\partial X}$.
By standard approximation theory (for example, see Hirsch~\cite{Hirsch}),
we may assume $F$ is smooth in the time direction and the fixed-time maps are smooth for time not equal to $0$.

By using a boundary collar in $\Sigma$ we can modify $h$ to get a map in $S_k$.
Let $N$ be a tubular neighborhood of $\partial\Sigma$ in $\Sigma$.
Let $\alpha\co N\to \partial\Sigma\times[0,1]$ be a diffeomorphism such that $\alpha(x)=(x,0)$ for $x\in\partial\Sigma$.
Let $\beta\co \overline{\Sigma-N}\to\Sigma$ be a diffeomorphism such that $\alpha(x)=(\beta(x),1)$ for $x\in\partial N-\partial\Sigma$.
Let $i$ be defined by:
\[i(x)=\left\{\begin{array}{ll}
h(\beta(x)) & \mbox{ if $x\in \overline{\Sigma-N}$}\\
F(\alpha(x)) & \mbox{ if $x\in N$}
\end{array}
\right.
\]
The reader can easily verify that $i\in S_k$.
\end{proof}

\begin{lemma}\label{le:Scontractible}
The space $S_k$ is contractible.
\end{lemma}

\begin{proof}
Let $\widetilde \Sigma$ be the universal cover of $\Sigma$ and pick a basepoint $*\in\widetilde\Sigma$ that projects to $*\in\Sigma$.
Let $\widetilde S_k$ be 
\[\widetilde S_k=\{\tilde\imath\co \widetilde\Sigma\to G_k|\tilde\imath(*)=e \mbox{ and $\tilde\imath$ is a lift of some $i\in S_k$}\}\]
with the compact--open topology.

Fix a map $\tilde\imath\in\widetilde S_k$.
Define $\widetilde F\co \widetilde S_k\times [0,1]\to \widetilde S_k$ by
\[\widetilde F(\tilde\jmath, t)(\tilde x)=\tilde\jmath(\tilde x)\big(\tilde\jmath(\tilde x)^{-1}\tilde\imath(\tilde x)\big)^t\]
where $\tilde\jmath\in \widetilde S_k$, $t\in[0,1]$, and $\tilde x\in\widetilde \Sigma$.
Here the exponent is given by the Lie group logarithm and exponential maps.

Let $p\co \pi\to\gam{k}$ denote the projection.  
Since the maps in $S_k$ induce $p$ on fundamental groups, we have: 
\[\tilde\jmath(\alpha\cdot \tilde x)=p(\alpha)\tilde\jmath(\tilde x)\]
for any $\tilde\jmath\in \widetilde S_k$, $\tilde x\in \widetilde \Sigma$, and $\alpha\in\pi$.

The homotopy $\widetilde F(\tilde\jmath,t)$ obeys the same rule:
\begin{align*}
\widetilde F(\tilde\jmath,t)(\alpha\cdot \tilde x)&= \tilde\jmath(\alpha \tilde x)(\tilde\jmath(\alpha\tilde x)^{-1}\tilde\imath(\alpha\tilde x))^t\\
&=p(\alpha) \tilde\jmath(\tilde x)(\tilde\jmath(\tilde x)^{-1}p(\alpha)^{-1}p(\alpha)\tilde\imath(\tilde x))^t\\
&=p(\alpha)\widetilde F(\tilde\jmath,t)(\tilde x)
\end{align*}

Since $\widetilde F(\tilde\jmath,t)$ is observes this rule, it descends to a continuous map
$F(j,t)\co  \Sigma\to X_k$.
Here $j\in S_k$ is the result of projecting $\tilde\jmath$ down to $S_k$.
This notation is justified since $F(j,0)=j$ and $F$ clearly varies continuously in $j$ and $t$.
Further, for each $j$, $(F(j,t))_*=p$ on fundamental groups.
One can easily see that $\widetilde F(\tilde\jmath, t)$ is constant on $\partial\widetilde\Sigma$ (with respect to $t$ and $\tilde\jmath$),
so for all $t$ and $j$, we have that $F(j,t)|_{\partial \Sigma}$ is indeed $i_
{\partial\Sigma}$.
So the map $F$ is actually a map
\[F\co S_k\times[0,1]\to S_k\]
and $F(j,1)=i$ for any $j$.
So $F$ is a contraction for $S_k$.
\end{proof}

%%%%%%%%%%%%%%%%%%%%%%%%%%%%%%%%%%%%%%%%%%%%%%%%%%%%%%%%%%%%%%%
\subsection{The main construction}\label{ss:CoreConstruction}
Now we are ready to give the definition of the crossed homomorphism extending the $k$th Johnson homomorphism.

At this point we fix a chain $T\in C_3(\Sigma\times[0,1])$ representing the fundamental class of $\Sigma\times[0,1]$ relative to $\partial(\Sigma\times[0,1])$. 
This fundamental class is the one corresponding to the orientation of $\Sigma\times[0,1]$ induced from the orientation on $\Sigma$.

The map $v\co C_*(X_k;\R)\to C_*(\lag{k})$ is as defined in Definition~\ref{de:vectorchainmap}.

\begin{definition}[The extension of the $k$th Morita homomorphism]
For $i\in S_k$,
the \emph{extended $k$th Morita map relative to $i$} is the crossed homomorphism
\begin{align*}
\emomap{k,i}:\bdM&\to C_3(\lag{k})/B_3(\lag{k})\\
\tag*{\text{defined by}}
[\phi]&\mapsto [v((F_\phi)_*T)]
\end{align*}
where $\phi$ is any representative of $[\phi]$ and the map $F_\phi\co \Sigma\times [0,1]\to X_k$ is any piecewise-smooth homotopy from $i$ to $\phi\cdot i$ through maps in $S_k$.

\end{definition}
The homotopy $F_\phi$ exists by Lemma~\ref{le:Scontractible}.
Although $\emomap{k,i}$ depends on the choice of $i\in S_k$, we will suppress this and write $\emomap{k}$ when there is no confusion.
We proceed by proving that the value of $\emomap{k}$ is independent of the other choices in its definition, and then by proving that it is, in fact, a crossed homomorphism.
The other choices in the definition are
 the choice of representative $T$ of $[\Sigma\times [0,1]]$,
 the choice of representative $\phi$ of $[\phi]$
 and the choice of homotopy $F_\phi$.
\begin{lemma}\label{le:2dimensional}
Suppose $M$ is a surface, and $f\co M\to X_k$ is piecewise-smooth.
Then for any chain $C\in C_3(M)$, we have $v(f_*C)=0$.
In particular, if $N$ is a $3$--manifold, $D\in C_3(N)$ and the piecewise-smooth map $F\co N\to X_k$ factors through $f$, then $v(F_*D)=0$. 
\end{lemma}

\begin{proof}
By the coordinate form in Lemma~\ref{le:CoordinateForm}, the chain $v(f_*C)$ is a linear combination of basis elements, where the coefficients are integrals of differential $3$--forms over $f_*C$.
Since $f$ is piecewise-smooth, we may pull back these differential $3$--forms to $M$ and compute the integrals there.
Of course, every $3$--form on $M$ is degenerate, so $v(f_*C)=0$.
The lemma follows.
\end{proof}

\begin{lemma}\label{le:lostcycle}
Suppose $F\co \Sigma\times[0,1]\to X_k$ is a piecewise-smooth homotopy through elements of $S_k$.
Then $[v(F_*T)]\in C_3(\lag{k})/B_3(\lag{k})$ does not depend on the choice $T$ of representative of the fundamental class of $\Sigma\times[0,1]$ relative to its boundary.
In particular, the value of $\emomap{k}([\phi])$ does not depend on the choice of representative of the fundamental class of $\Sigma\times [0,1]$ relative to its boundary.
\end{lemma}

\begin{proof}
Let $T$ and $T'$ represent the fundamental class of $\Sigma\times[0,1]$ relative to its boundary.
Then $\partial T$ and $\partial T'$ represent the same class in $H_2(\partial(\Sigma\times[0,1]))$.
Pick a class $C\in C_3(\partial(\Sigma\times[0,1]))$ with $\partial C = \partial T-\partial T'$.
Then the chain $T-T'-C$ is a cycle in $Z_3(\Sigma\times[0,1])$.
Furthermore, since $H_3(\Sigma\times[0,1])=H_3(\Sigma)=0$, we know $T-T'-C\in B_3(\Sigma\times[0,1])$.

The chain $C$ is supported on $\partial(\Sigma\times [0,1])$, 
so by Lemma~\ref{le:2dimensional}:
\[v(F_*C)=0\]
The claim follows from:
\[v(F_*T)-v(F_*T')= v(F_*(T-T'-C))\in v(B_3(X_k))\subset B_3(\lag{k})
%\proved
\]
\end{proof}

\begin{lemma}\label{le:lostdiffeo}
The value of $\emomap{k}([\phi])$ does not depend on the choice of homotopy $F_\phi$ in its definition, for a fixed choice of representative $\phi$ of $[\phi]$.
\end{lemma}
\begin{proof}
Let $i'=\phi\cdot i$.
Suppose that $F,F'\co \Sigma\times [0,1]\to X_k$ are two homotopies through maps in $S_k$, both from $i$ to $i'$.
It follows from Lemma~\ref{le:Scontractible} that there is a homotopy 
$K\co \Sigma\times[0,1]^2\to X_k$ through maps in $S_k$, from $F$ to $F'$ relative to their endpoints.

We can pick a representative $C\in C_4(\Sigma\times[0,1]^2)$ of the fundamental class of $\Sigma\times[0,1]^2$ relative to its boundary such that:
\[\partial K_* C = F'_*T+\underline{i'}_*T-F_*T-\underline{i}_*T+(\underline{i_{\partial\Sigma}})_*(\partial T\times[0,1]^2)\]
Here $\underline{i}$, $\underline{i'}$ and $\underline{i_{\partial\Sigma}}$ denote the constant homotopies and $\partial T\times [0,1]^2$ is a representative of the fundamental class of $\partial\Sigma\times[0,1]^2$ relative to boundary.
By Lemma~\ref{le:2dimensional}:
\[v(\underline{i'}_*T)=v(\underline{i}_*T)=v((\underline{i_{\partial\Sigma}})_*(\partial T\times[0,1]^2))=0\]
Then we have
\[v(F'_*T)-v(F_*T)=v(\partial K_*C)=\partial v(K_*C)\]
which is in $B_3(\lag{k})$.
\end{proof}

\begin{lemma}\label{le:crossedhom}
Let $F$ and $F'$ be piecewise-smooth homotopies through maps in $S_k$ such that $F_1=F'_0$.
Let $F''$ be the concatenation of $F$ with $F'$.
Then $F''$ is a piecewise-smooth homotopy through maps in $S_k$ and we have:
\[[v(F_*T)]+[v(F'_*T)]=[v(F''_*T)]\]
\end{lemma}

\begin{proof}
It is obvious that $F''$ is a piecewise-smooth homotopy through maps in $S_k$.

Choose a representative $T'\in C_3(\Sigma\times [0,1])$ of the fundamental class relative to boundary, such that:
\[F''_*T'=F_*T+F'_*T\]
By Lemma~\ref{le:lostcycle}, we have:
\[[v(F''_*T)]=[v(F''_*T')]\]
The lemma follows.
\end{proof}

\begin{lemma}\label{le:constructionequivariant}
Suppose that $F$ is a piecewise-smooth homotopy through maps in $S_k$ and $\phi\in\Diff(\Sigma,\partial\Sigma)$.
Then:
\[[v((\phi\cdot F)_*T)]=\phi\cdot [v(F_*T)]\]
\end{lemma}

\begin{proof}
By the definitions of the actions, we have:
\[(\phi\cdot F)_*T=\phi\cdot (F_*\phi^{-1}_*T)\]
However, since $\phi$ is orientation preserving, $\phi^{-1}_*T$ is also a representative of the fundamental class of $\Sigma\times[0,1]$ relative to the boundary.
The lemma follows from Lemma~\ref{le:lostcycle}.
\end{proof}
\begin{claim}
The value of $\emomap{k}([\phi])$ does not depend on the choice of representative $\phi$ of $[\phi]$.
\end{claim}

\begin{proof}
Suppose $\phi,\psi\in\Diff(\Sigma,\partial\Sigma)$ are isotopic.
Then there is a smooth homotopy $K$ on $\Sigma$, relative to $*$, taking the identity map to $\phi\psi^{-1}$.
Then $i\circ K$ is a smooth homotopy from $i$ to $\phi\psi^{-1}\cdot i$ through maps in $S_k$.
Lemma~\ref{le:2dimensional} tells us $v((i\circ K)_*T)=0$.

Let $F_\psi$ be a piecewise-smooth homotopy from $i$ to $\psi\cdot i$ through maps in $S_k$, and let $F'$ be the concatenation of $i\circ K$ with $\phi\psi^{-1}\cdot F_\psi$.
Then $F'$ is a piecewise-smooth homotopy from $i$ to $\phi\cdot i$, and by Lemma~\ref{le:crossedhom} and Lemma~\ref{le:constructionequivariant} we have:
\[[v(F'_*T)]=[v((i\circ K)_*T)]+[v((F_\psi)_*T)]\]
Since $v((i\circ K)_*T)=0$, we get the same value for $\emomap{k}([\phi])$ using either representative.
\end{proof}

\begin{claim}
The map $\emomap{k}$ is a crossed homomorphism.
\end{claim}

\begin{proof}
Let $\phi,\psi\in\Diff(\Sigma,\partial\Sigma)$.
Let $F_\phi, F_\psi$ be piecewise-smooth homotopies through maps in $S_k$ from $i$ to $\phi\cdot i$ and  $\psi\cdot i$ respectively.
Let $F'$ be the concatenation of $F_\phi$ with $\phi\cdot F_\psi$.
Then $F'$ is a piecewise-smooth homotopy through maps in $S_k$ from $i$ to $\phi\psi\cdot i$.
By Lemma~\ref{le:crossedhom} and Lemma~\ref{le:constructionequivariant}, we have:
\[[v(F'_*T)]=[v((F_\phi)_*T)]+\phi\cdot [v((F_\psi)_*T)]\]
The claim follows immediately.
\end{proof}

The previous claim completes the proof that $\emomap{k}$ is well-defined.

\begin{proposition}\label{pr:lostbasemap}
The cohomology class \[[\emomap{k}]=[\emomap{k,i}]\in H^1(\bdM,C_3(\lag{k})/B_3(\lag{k}))\]
does not depend on the choice of $i\in S_k$, or on the choice of $i_{\partial\Sigma}$ in the definition of $S_k$.
\end{proposition}

\begin{proof}
First we show the independence from the choice of $i\in S_k$.
Suppose $i,i'\in S_k$.
Let $K$ be any smooth homotopy from $i$ to $i'$ through maps in $S_k$.
For $\phi\in\Diff(\Sigma,\partial\Sigma)$, let $F'_\phi$ be a piecewise-smooth homotopy from $i'$ to $\phi\cdot i'$,  through maps in $S_k$.
Let $F_K$ be the concatenation of $K$ with $F'_\phi$ and then with the homotopy $\phi\cdot K$ travelled in reverse.

This $F_K$ is a piecewise-smooth homotopy from $i$ to $\phi\cdot i$ through maps in $S_k$.
By the definition of $\emomap{k}$, and Lemma~\ref{le:crossedhom} and Lemma~\ref{le:constructionequivariant}:
\begin{align*}
\emomap{k,i}([\phi])&=[v((F_K)_*T)]\\
&=[v(K_*T)]+[v((F'_\phi)*T)]-[v((\phi\cdot K)_*T)]\\
&=\emomap{k,i'}([\phi])-([\phi]\cdot [v(K_*T)]-[v(K_*T)])\\
&=\emomap{k,i'}([\phi])-(\delta[v(K_*T)])([\phi])
\end{align*}
Since $\delta[v(K_*T)]\in B^1(\bdM,C_3(\lag{k})/B_3(\lag{k}))$, we have shown that $[\emomap{k}]$ does not depend on $i$ if $i_{\partial\Sigma}$ is fixed.

Now suppose that $i_{\partial\Sigma}$ and $i_{\partial\Sigma}'\co \partial\Sigma\to \ell_k$ are two different maps as in the definition of $S_k$, and suppose that $S_k$ and $S_k'$ are the spaces of maps they respectively define.
Let $i\in S_k$ and 
let $K\co \partial\Sigma\times[0,1]\to \ell_k$ be a homotopy relative to $*$ from $i_{\partial\Sigma}$ to $i_{\partial\Sigma}'$.
We can get an element $i'\in S_k'$ by extending $i$ by $K$ on a collar of the boundary of $\Sigma$.
This is the same trick as in the proof of Lemma~\ref{le:Snonempty}, so we omit the details.

For $\phi\in\Diff(\Sigma,\partial\Sigma)$, let $F_\phi$ be a piecewise-smooth homotopy through maps in $S_k$, with $F_\phi$ going from $i$ to $\phi\cdot i$.
Let $\underline K\co \partial\Sigma\times[0,1]^2\to\ell_k$ be the constant homotopy from $K$ to itself.
Let $N$ be a tubular neighborhood of $\partial \Sigma$. 
Let $F'_\phi$ be a piecewise-smooth homotopy from $i'$ to $\phi\cdot i'$ through maps in $S_k$, such that $F'_\phi|_{N\times[0,1]}$ factors through $\underline K$ and $F'_\phi|_{\overline{\Sigma- N}\times[0,1]}$ factors through $F_\phi$.
Let $T'\in C_3(\Sigma\times [0,1])$ and $C\in C_3(\partial\Sigma\times[0,1]^2)$ be representatives for the respective fundamental classes relative to boundary such that:
\[(F'_\phi)_*T'=(F_\phi)_*T+\underline K_*C\]
By Lemma~\ref{le:2dimensional}, we know $v(\underline K_* C)=0$.
Then we get the same value for $\emomap{k}([\phi])$ using $i_{\partial\Sigma}$ and $i$ as we get using $i_{\partial\Sigma}'$ and $i'$.
The proposition follows.
\end{proof}

\subsection{The main construction restricted to the $k$th Torelli group}
\label{ss:topologicalmoritahomomorphism}
The goal of this section is to show that the main construction works nicely with simplicial homology when we restrict to the $k$th Torelli group $\bdI{k}$.

Let $C_\Sigma$ be a representative of the fundamental class of $\Sigma$ relative to $\partial\Sigma$.
By abuse of notation, let $C_\Sigma\times\{0\}$ and $C_\Sigma\times\{1\}$ represent the images of $C_\Sigma$ under the identifications of $\Sigma$ with $\Sigma\times\{0\}$ and $\Sigma\times\{1\}$.
Similarly, let $\partial C_\Sigma\times[0,1]$ and $C_\Sigma\times[0,1]$ be representatives of the fundamental classes of $\partial\Sigma\times[0,1]$ and $\Sigma\times[0,1]$ relative to boundaries, respectively, so that we have:
\[\partial(C_\Sigma\times[0,1])=\partial C_\Sigma\times[0,1]+C_\Sigma\times\{1\}-C_\Sigma\times\{0\}\]
In this section $H_*(X_k)$ denotes the piecewise-smooth simplicial homology with integer coefficients.

\begin{definition}\label{de:topologicalversion}
The \emph{topological version of the $k$th Morita homomorphism} is the map
\begin{align*}
\tomohom{k}\co \bdI{k}&\to H_3(X_k)\\
\tag*{\text{given by:}} 
[\phi]&\mapsto [(F_\phi)_*(C_\Sigma\times[0,1])-i_*D_\phi]
\end{align*}
Here 
$i\in S_k$,
$\phi$ is a representative of $[\phi]$,
the map $F_\phi\co \Sigma\times [0,1]\to X_k$ is a piecewise-smooth homotopy from $i$ to $\phi\cdot i$ through maps in $S_k$ and
the chain $D_\phi\in C_3(\Sigma)$ bounds the chain 
\[\phi^{-1}_*C_\Sigma-C_\Sigma+(proj_1)_*(\partial C_\Sigma\times [0,1])\]
where $proj_1\co \partial\Sigma\times [0,1]\to \partial\Sigma$ is the first coordinate projection. 
\end{definition}

\begin{claim}
The chain $(F_\phi)_*(C_\Sigma\times[0,1])-i_*D_\phi$ in Definition~\ref{de:topologicalversion} exists and is a cycle.
\end{claim}

\proof
The chain $D_\phi$ exists because its specified boundary is a cycle in $Z_2(\Sigma)$ and $H_2(\Sigma)=0$.
The claim then follows from this computation:
\begin{align*}
\partial ((F_\phi)_*C_\Sigma\times[0,1])&=(F_\phi)_*(\partial C_\Sigma\times[0,1]+C_\Sigma\times\{1\}-C_\Sigma\times\{0\})\\
&=i_*(proj_1)_*(\partial C_\Sigma\times[0,1])+i_*\phi^{-1}_*C_\Sigma-i_*C_\Sigma\\
&= \partial i_*D_\phi
\rlap{\hspace{2.8 in}\qedsymbol }
\end{align*}

\begin{lemma}\label{le:topologytoalgebra}
For any $i\in S_k$, the maps $v_*\circ \tomohom{k}$ and $\emomap{k,i}|_{\bdI{k}}$ are equal.
\end{lemma}

\begin{proof}
Fix $\phi\in\Diff(\Sigma,\partial\Sigma)$ with $[\phi]\in\bdI{k}$.
We have $v(i_*D_\phi)=0$ by Lemma~\ref{le:2dimensional}.
Then 
\[v_*(\tomohom{k,i}([\phi]))=[v((F_\phi)_*(C_\Sigma\times[0,1]))]=\emomap{k,i}([\phi])\]
where $F_\phi$ is as in definition of $\tomohom{k,i}$.
\end{proof}

\begin{lemma}\label{le:homologyinjection}
The map $H_3(X_k)\to H_3(X_k;\R)$ is an injection.
\end{lemma}

\begin{proof}
A theorem of Igusa and Orr~\cite[Theorem~5.9]{Igusa&Orr} states that $H_3(\gam{k})$ is torsion free.
The claim follows.
\end{proof}

\begin{claim}
The map $\tomohom{k}$ is a well-defined homomorphism and does not depend on the choices in its definition, including the choices of $i\in S_k$ and $i_{\partial\Sigma}$ in the definition of $S_k$.
\end{claim}

\begin{proof}
By Lemma~\ref{le:homologyinjection} and Corollary~\ref{co:dualofNomizuTheorem}, the map $v_*\co H_3(X_k)\to H_3(\lag{k})$ is an injection.
Then by Lemma~\ref{le:topologytoalgebra}, and the fact that $\emomap{k}$ is well-defined, we see that $\tomohom{k}$ does not depend on any of the choices except possibly $i$ and $i_{\partial\Sigma}$.
The fact that $v$ is equivariant implies that $\tomohom{k}$ is a crossed homomorphism, but the action of $\bdI{k}$ on $\CmodB{k}$ is trivial, making $\tomohom{k}$ a homomorphism.
By Proposition~\ref{pr:lostbasemap}, varying $i$ or $i_{\partial\Sigma}$ only changes $\tomohom{k}$ by adding a coboundary, but $B^1(\bdI{k},\CmodB{k})$ is trivial.
\end{proof}

\begin{remark}
It is a difficult fact that third homology of $X_k$ is isomorphic to the $3$--dimensional relative (non-spin) bordism group $\Omega_3(\gam{k})$.
In light of this fact, this $\tomohom{k}$ can be identified with a homomorphism on $\bdI{k}$ defined by Heap in~\cite{Heap}.
We will not pursue this connection in the current paper. 
\end{remark}
%%%%%%%%%%%%%%%%%%%%%%%%%%%%%%%%%%%%%%%%%%%%%%%%%%%%%%%%%%%%%%%%%%%%%%%%%%
\section{Equivalence of $\tomohom{k}$ and Morita's homomorphism}\label{se:HomomorphismEquivalence}
\subsection{Morita's homomorphism}\label{ss:Moritahomomorphism}

We recall some ideas from group homology theory, for which Brown~\cite{Brown} is a reference.
Let $G$ be a group and let $F_*(G)$ and $C_*(G)$ denote the standard (bar) resolution and standard group chains, respectively.
Of course, $H_*(G)=H_*(C_*(G))$.
The $G$--complex $F_*(G)$ has a $\Z$--basis of elements of the form $(g_0,\ldots,g_n)$.
This descends to a basis for $C_*(G)$ of elements of the form $[(g_0,\ldots,g_n)]$, where the brackets denote the co-invariance class.
A homomorphism $f$ from $G$ to a group $G'$ induces a map $f_*\co F_*(G)\to F_*(G')$, given by  $[(\gamma_0,\ldots,\gamma_n)]\mapsto[(f(\gamma_0),\ldots,f(\gamma_n))]$.
This gives the action of $\Aut(G)$ on $F_*(C)$,
which induces the action on $C_*(G)$ and in turn induces the usual action on $H_*(G)$.

Now consider the complex $C_*(\pi)$.  
Since $\bdM$ injects in $\Aut(\pi)$, we have $\bdM$ acting on $C_*(\pi)$.
Let $\ell\in\pi$ be the class of the positive boundary loop in $\Sigma$.
Then clearly, $[\ell]\in B_1(\pi)$ is fixed by the action of $\bdM$.
Let $p\co \pi\to\gam{k}$ be the projection.
The map $p$ and the induced map $p_*\co C_*(\pi)\to C_*(\gam{k})$ are both obviously $\bdM$--equivariant.

As mentioned in the introduction, we have a homomorphism called the \emph{$k$th Morita homomorphism}.
To reiterate, Morita's homomorphism 
\[\mohom{k}\co \bdI{k}\to H_3(\gam{k})\]
is defined by:
\[[\phi]\mapsto [p_* D_{[\phi]}]\]
Here the chain $D_{[\phi]}\in C_3(\pi)$ bounds $[\phi]\cdot C- C$,
where $C\in C_2(\pi)$ is a chain bounding the element $-[\ell]\in B_1(\pi)$ and projecting to a representative the fundamental class of $\pi/\overline{\langle\ell\rangle}$.

It makes sense to talk about the homology class of $p_*D_{[\phi]}$ when $[\phi]\in\bdI{k}$;
 because $\bdI{k}$ acts trivially on $\gam{k}$ we have:
\[\partial p_* D_{[\phi]} = p_*\partial D_{[\phi]}=p_*[\phi]\cdot C - p_*C=[\phi]\cdot p_*C-p_*C=0\]
For proofs that this definition does not depend on the choices involved and that it is a homomorphism, we refer the reader to Morita's paper~\cite[Theorem~3.1]{Morita93}.

\subsection{Equivalence of the homomorphisms}\label{ss:HomomorphismEquivalence}
In this subsection we show:
\begin{theorem}\label{th:homomorphismequivalence}
For each $k\geq 2$, 
\[\mohom{k}=f\circ \tomohom{k}\]
where \mbox{$f\co H_3(X_k)\to H_3(\gam{k})$} is the canonical isomorphism.
\end{theorem}

In order to prove this, we find a chain equivalence between $C_*(X_k)$ and $C_*(\gam{k})$
and show that for each $[\phi]\in\bdI{k}$, this equivalence maps a representative of $\tomohom{k}([\phi])$ in $C_3(X_k)$ to a representative of $\mohom{k}([\phi])$ in $C_3(\gam{k})$.  

We recall that a singular $n$--simplex $\Delta\in C_n(X)$ on a space $X$ is a map from a fixed $n$--simplex (considered as a space) to $X$.
We denote the vertices of our fixed $n$--simplex as $v_0,\ldots, v_n$, so that the vertices of the $n$--simplex $\Delta$ are the points $\Delta(v_0),\ldots,\Delta(v_n)\in X$.

\begin{definition}\label{de:chainequivalence}
Suppose the group $G$ is the fundamental group of an aspherical space $X$, and $\widetilde X\to X$ the universal cover of $X$.
Suppose $A\subset \widetilde X$ is a fundamental domain for the $G$ action on $\widetilde X$.

Let $c_A\co \widetilde X\to G$ to be the map sending $\tilde x$ to the unique $g\in G$ with $\tilde x\in g\cdot A$.
Let $\tilde f_A\co C_*(\widetilde X)\to F_*(G)$ be the $G$--complex map given by
\[\tilde f_A(\Delta)=(c_A(\Delta(v_0)),\ldots,c_A(\Delta(v_n))),\]
for $\Delta$ an $n$--simplex in $\widetilde X$.
Let $f_A\co C_*(X)\to C_*(G)$ be the chain complex map given by
\[f_A(C)=[\tilde f_A(\widetilde C)],\] 
for $C\in C_*(X)$, where  $\widetilde C$ is a lift of $C$ to $\widetilde X$.
\end{definition}

\begin{lemma}\label{le:reallyequivalence}
The maps $\tilde f_A$ and $f_A$ are well-defined. 
Furthermore, both maps are homotopy equivalences.
\end{lemma}

\begin{proof}
We leave it to the reader to check $\tilde f_A$ is a $G$--complex map and that $f_A$ is a chain complex map.
The map $f_A$ does not depend on the choice of $\widetilde C$ because the $G$--action is lost in passing to $C_*(G)$.

The set of $n$--simplices $\Delta\in C_*(\widetilde X)$ such that $\Delta(v_0)\in A$ form a free $\Z G$--basis of $C_n(\widetilde{X})$.
Since $\widetilde X$ is contractible, the $G$--complex $C_*(\widetilde X)$ is a free resolution of $\Z$.
Since $C_*(\widetilde{X})$ and $F_*(G)$ are both free resolutions of $\Z$ and $\tilde f_A$ preserves their respective augmentation maps $C_0(\widetilde X)\to \Z$ and $C_0(G)\to \Z$, it follows from standard theory (see Brown~\cite[chapter I.7]{Brown}) that $\tilde f_A$ is a homotopy equivalence of $G$--complexes.  
It immediately follows that $f_A$ is a homotopy equivalence of chain complexes.
\end{proof}

\begin{lemma}\label{le:constantimage}
Given the setup of Definition~\ref{de:chainequivalence}, suppose that $\tilde x\in A$ is a basepoint for $\widetilde X$.
Let $x\in X$ be the image of $\tilde x$ under the covering projection.
Suppose $\Delta_t$ for $t\in [0,1]$ is a homotopy of a simplex of $C_*(X)$ in $X$,
subject to the condition that 
\[\Delta_t(v_i)=x\]
for each vertex $v_i$, for $0\leq t\leq 1$.
Then $f_A(\Delta_t)$ is constant with respect to $t$.
\end{lemma}

\begin{proof}
For any $i\in\{1,\ldots, n\}$, we have an element $\gamma_{i}\in \pi_1(X,x)=G$ equal to the class of the edge of $\Delta_t$ from $v_0$ to $v_i$, considered as a loop in $X$.
Since $\Delta_t$ and $\Delta_{t'}$ are homotopic relative to the vertices for all $t$ and $t'$, these elements $\gamma_{i}$ do not depend on $t$.
For each $t$, let $\widetilde \Delta_t$ be the lift of $\Delta_t$ to $\widetilde X$ with the $\widetilde \Delta_t(v_0)=\tilde x$.
Then $\widetilde \Delta_t(v_i)=\gamma_i\cdot \tilde x$ and $c_A(\widetilde \Delta_t(v_i))=\gamma_i$ for all $t$.
The lemma follows from the definition of $f_A$.
\end{proof}

\begin{lemma}\label{le:commutingdiagram}
Suppose $i\in S_k$.
Let $\widetilde\Sigma$ be the universal cover of $\Sigma$.
It is possible to find a fundamental domain $A$ for the action of $\pi$ on $\widetilde \Sigma$, and a fundamental domain $B$ for the action of $\gam{k}$ on $G_k$, 
such that the following diagram commutes:
\[
\xymatrix{
C_*(\Sigma)\ar[r]^{f_A}\ar[d]^{i_*}&C_*(\pi)\ar[d]^{p_*}\\
C_*(X_k)\ar[r]^{f_B}&C_*(\gam{k})
}
\]
\end{lemma}

\begin{proof}
Let $\tilde\imath\co \widetilde \Sigma\to G_k$ be a lift of $i$.
If $\gamma\in\pi$ and $\tilde x\in\widetilde \Sigma$, we have:
\[\tilde\imath(\gamma\cdot \tilde x)=p(\gamma)\tilde\imath(\tilde x)\]
This follows from the definition of $S_k$.
Since $\tilde\imath$ observes this rule, it descends to a map $\hat\imath\co \widetilde \Sigma/\ker p\to G_k$.
Let $B\subset G_k$ be a fundamental domain for the action of $\gam{k}$ on $G_k$.
Then $\hat\imath^{-1}(B)$ is a fundamental domain for the action of $\gam{k}$ on $\widetilde \Sigma/\ker p$.
Let $A'\subset \widetilde\Sigma$ be a fundamental domain for the action of $\ker p$ on $\widetilde \Sigma$.
Set $A=A'\cap \tilde\imath^{-1}(B)$.

We claim that $A$ is a fundamental domain for the action of $\pi$ on $\widetilde \Sigma$.
If $\tilde x\in \widetilde \Sigma$, then there is a unique $\gamma\in\gam{k}$ with $\gamma\tilde\imath(\tilde x)\in B$.
If $\tilde\gamma\in\pi$ is a lift of $\gamma$, then there is a unique $\alpha\in\ker p$ with $\alpha\tilde\gamma\tilde x\in A$.
It follows that $\alpha\tilde\gamma$ is the unique element of $\pi$ sending $\tilde x$ into $A$.
Since $\tilde x$ was arbitrary, this means that $A$ is a fundamental domain.

Now we show that the diagram commutes.
Pick $\tilde x\in \widetilde \Sigma$.
Let $\alpha= c_A(\tilde x)$.
Then $\tilde\imath(\tilde x)\in \tilde\imath(\alpha\cdot A)$.
We deduce that $\tilde\imath(\tilde x)\in p(\alpha) B$, so $c_B(\tilde\imath(\tilde x))=p(\alpha)$.
We have shown: 
\[c_B\circ\tilde\imath = p\circ c_A\]
It follows from Definition~\ref{de:chainequivalence} that the diagram commutes.
\end{proof}

At this point, we fix a representative $\ell\in Z_1(\partial\Sigma)\subset Z_1(\Sigma)$ of the boundary loop in $\Sigma$ and a class $C_{\Sigma}\in C_2(\Sigma)$ with $\partial C_{\Sigma}=\ell$.
We demand that the simplices that are summands of $C_{\Sigma}$ have all of their vertices equal to the basepoint $*$ of $\Sigma$.
The symbols $C_\Sigma\times[0,1]$ 
and $\partial C_\Sigma\times[0,1]$ have the same (abusive) meanings that they were given in the preamble to Definition~\ref{de:topologicalversion}.
We also fix a map $i_{\partial\Sigma}\co \partial\Sigma\to X_k$, which determines the space $S_k$, and we fix a map $i\in S_k$ as in Definition~\ref{de:topologicalversion}.

\begin{lemma}\label{le:losthomotopy}
Suppose $[\phi]\in \bdI{k}$ and $F$ is the chain 
\[F=(F_\phi)_*(C_{\Sigma}\times[0,1])\in C_3(X_k)\]
where $F_\phi\co \Sigma\times[0,1]\to X_k$ is a piecewise-smooth homotopy from $i$ to $\phi\cdot i$ through maps in $S_k$ as in Definition~\ref{de:topologicalversion}.
Suppose $B$ is a fundamental domain for the action of $\gam{k}$ on $X_k$.
Then $f_B(F)$ is a boundary.
\end{lemma}

\begin{proof}
By definition, the chain $C_{\Sigma}$ is a sum of simplices whose vertices are all equal to the basepoint $*$ of $\Sigma$.
Since the homotopy $F_{\phi}$ is a homotopy relative to the basepoint $*$, the homotopy $F_{\phi,t}(\Delta)$ satisfies the conditions of Lemma~\ref{le:constantimage}, where $\Delta$ varies over the summands of $C_{\Sigma}$.
This means that $f_B(F_{\phi,t}(C_{\Sigma}))$ is constant with respect to $t$.
So since $f_B(F)$ is the image of a constant homotopy, it is a boundary.
\end{proof}

\begin{lemma}\label{le:equivariantequivalence}
Let $A$ be a fundamental domain for the action of $\pi$ on $\widetilde \Sigma$.
If $\phi$ is in $\Diff(\Sigma,\partial\Sigma)$ and $[\phi]$ is its image in $\bdM$, we have:
\[f_A(\phi\cdot C_{\Sigma})=[\phi]\cdot f_A(C_{\Sigma})\]
\end{lemma}

\begin{proof}
Again, we use the fact that the summands of $C_{\Sigma}$ have all vertices equal to the basepoint $*$ of $\Sigma$.
We restrict our attention to a single simplex $\Delta$ of $C_{\Sigma}$.
We note that, as in the proof of Lemma~\ref{le:constantimage}, the image of this simplex under $f_A$ is the co-invariance class of the $n$--tuple of elements of $\pi$ specified by the edges of $\Delta$ (considered as loops) from $v_0$ to $v_i$ for $i=0,\ldots,n$.
Since the action of $\phi$ on loops in $\Sigma$ defines the action of $[\phi]$ on $\pi$, this proves the lemma.
\end{proof}

\begin{proof}[Proof of Theorem~\ref{th:homomorphismequivalence}]
Let $A$ and $B$ be fundamental domains for the actions of $\pi$ on $\widetilde \Sigma$ and $\gam{k}$ on $X_k$ respectively, such that the conclusions of Lemma~\ref{le:commutingdiagram} hold.
Let $\phi\in\Diff(\Sigma,\partial\Sigma)$ be an element such that $[\phi]\in\bdI{k}$.
Let $D_\phi\in C_3(\Sigma)$ be some chain satisfying
\[\partial D_\phi=\phi^{-1}_*C_\Sigma-C_\Sigma+(proj_1)_*(\partial C_\Sigma\times [0,1])\]
where $proj_1\co \partial\Sigma\times [0,1]\to \partial\Sigma$ is the first coordinate projection. 
Let $F_\phi\co \Sigma\times[0,1]\to X_k$ be a piecewise-smooth homotopy from $i$ to $\phi\cdot i$ through $S_k$, and let
\[C_\phi=(F_\phi)_*(C_\Sigma\times[0,1])-i_*D_\phi\]
so that
$\tomohom{k}([\phi])=[C_\phi]$.

%Apply Lemma~\ref{le:losthomotopy} and Lemma~\ref{le:commutingdiagram} to compute $f_B(C_\phi)$:
By Lemma~\ref{le:commutingdiagram}, we have:
\begin{align*}
f_B(C_\phi-(F_\phi)_*(C_{\Sigma}\times[0,1]))&=-f_B(i_* D_\phi)\\
&=-p_*f_A(D_\phi)
\end{align*}

Lemma~\ref{le:equivariantequivalence} allows us to compute $\partial f_A(D_\phi)$:
\begin{align*}
\partial f_A(D_\phi)&=f_A(\partial D_\phi)\\
&=f_A(\phi\cdot C_{\Sigma}-C_{\Sigma})\\
&=[\phi]\cdot f_A(C_{\Sigma})-f_A(C_{\Sigma})
\end{align*}

Then $f_B(C_\phi-(F_\phi)_*(C_{\Sigma}\times[0,1]))$ is $p_*(-[\phi]\cdot f_A(C_{\Sigma})+f_A(C_{\Sigma}))$ where $-f_A(C_{\Sigma})$ bounds the reverse boundary loop element $-f_A(\ell)$.
This fits the definition of $\mohom{k}$, so we have shown:
\[\mohom{k}([\phi])=[f_B(C_\phi-(F_\phi)_*(C_{\Sigma}\times[0,1]))]\]
Then by Lemma~\ref{le:losthomotopy}, we are finished:
\[\mohom{k}([\phi])=[f_B(C_\phi)]
%\proved
\]
\end{proof}

\begin{proof}[Final step in the proof of Main Theorem~\ref{mt:ExtMoritaHom}]
By Lemma~\ref{le:homologyinjection}, we know that the map $v_*\circ f_{B*}^{-1}\co  H_3(\gam{k})\into \CmodB{k}$ is indeed an injection.
By Lemma~\ref{le:topologytoalgebra}, we know that $\emomap{k}|_{\bdI{k}}=v_*\circ \tomohom{k}$ and 
by Theorem~\ref{th:homomorphismequivalence}, we know that $f_{B*}^{-1}\circ \mohom{k}=\tomohom{k}$.
This proves the theorem. 
\end{proof}

%%%%%%%%%%%%%%%%%%%%%%%%%%%5
\section{Extending the higher Johnson homomorphisms}
\label{se:proofof}
Recall from the introduction that $\lv{k+1}=\pi^{(k-1)}/\pi^{(k)}$.  We have:
\begin{equation}1\to\lv{k+1}\to\gam{k+1}\to\gam{k}\to 1\label{eq:groupsequence}\end{equation}
By taking Mal'cev completions (see Theorem~\ref{th:CompleteMaps}) we have an exact sequence of Lie groups:
\[1\to\lvr{k+1}\to G_{k+1}\to G_k \to 1\]
Taking $\lal{k+1}$ to be the Lie algebra of $\lvr{k+1}$, we get the following exact sequence of Lie algebras:
\begin{equation}0\to\lal{k+1}\to\lag{k+1}\to\lag{k}\to 0\label{eq:algebrasequence}\end{equation}
The Hochschild--Serre spectral sequence for the homology of this exact sequence has a differential:
\begin{align*}
d^2\co H_3(\lag{k})\to H_1(\lag{k},H_1(\lal{k+1}))&\cong \lag{2}\otimes\lal{k+1}\\
&\cong (\lag{k+1}\wedge\lal{k+1})/(\lag{k+1}^{(1)}\wedge\lal{k+1})
\end{align*}
Here $\lag{k+1}^{(1)}$ is the commutator subalgebra of $\lag{k+1}$, and $\lag{2}=H\otimes \R$ is the Lie algebra of $G_2=H\otimes \R$.
In short, one gets $d^2$ of a class $[c]$ in $H_3(\lag{k})$ by picking a representative $c$ in $Z_3(\lag{k})$, lifting $c$ to $\tilde c$ in $C_3(\lag{k+1})$, and taking the boundary of the lift.
This boundary is in $\lag{k+1}\wedge\lal{k+1}$ and is well-defined up to elements of $\lag{k+1}^{(1)}\wedge\lal{k+1}$.

We start by extending this differential to a map on $\CmodB{k}$.
\begin{definition}
The \emph{extended differential} 
\[\tilde{d^2}\co C_3(\lag{k})/B_3(\lag{k})\to C_2(\lag{k+1})/(\lag{k+1}^{(1)}\wedge \lal{k+1})\]
is defined by
\[\tilde{d^2}([c])=\partial\tilde c+\lag{k+1}^{(1)}\wedge\lal{k+1}\]
where $[c]\in C_3(\lag{k})/B_3(\lag{k})$, the chain $c\in C_3(\lag{k})$ is a representative of $[c]$, and the chain $\tilde c\in C_3(\lag{k+1})$ is a lift of $c$.
\end{definition}

\begin{claim}
The extended differential $\tilde{d^2}$ is well defined and extends the differential
$d^2\co H_3(\lag{k})\to (\lag{k+1}\wedge\lal{k+1})/(\lag{k+1}^{(1)}\wedge\lal{k+1}).$
\end{claim}

\begin{proof}
The fact that $\tilde{d^2}$ extends $d^2$ is immediate from the definitions.

To see that $\tilde{d^2}$ is well-defined, we first show that it does not depend on the choice of lift.
Suppose $c\in C_3(\lag{k})$ and that $\tilde c$ and $\tilde c'$ are lifts of $c$ to $C_3(\lag{k+1})$.
Then $\tilde c-\tilde c'$ is in the $\ker(C_3(\lag{k+1})\to C_3(\lag{k}))$.
This kernel is easily seen to be $\lag{k+1}\wedge\lag{k+1}\wedge\lal{k+1}$.
The formula for the boundary (Equation~(\ref{eq:KoszulBoundary})) and the fact that $\lal{k+1}$ is central in $\lag{k+1}$ indicate that image of this kernel under the boundary is in $\lag{k+1}^{(1)}\wedge\lal{k+1}$.
Therefore the choice of lift is irrelevant.

To show that the choice of representative of the class $[c]\in C_3(\lag{k})/B_3(\lag{k})$ is irrelevant, we show that the group of boundaries maps to zero.
Take $b\in B_3(\lag{k})$.
Pick a class $c\in C_4(\lag{k})$ with $\partial c=b$.
Pick a lift $\tilde c\in C_4(\lag{k+1})$.
Then $\tilde b = \partial \tilde c\in C_3(\lag{k})$ is a lift of $b$.
Immediately, we get that $\partial\tilde b=\partial\partial \tilde c = 0$, which proves the claim.
\end{proof}

The subspace $\lag{k+1}^{(1)}\wedge\lal{k}$ is clearly invariant under the action of $\bdM$ on $C_2(\lag{k})$.
This gives us an induced action of $\bdM$ on $C_2(\lag{k+1})/(\lag{k+1}^{(1)}\wedge \lal{k+1})$.
\begin{lemma}\label{le:d2equivariant}
The extended differential $\tilde{d^2}$ is $\bdM$--equivariant.
\end{lemma}
\begin{proof}This is immediate from the definition.
\end{proof}

Let $H^*$ denote $\Hom(H,\Z)$, the dual to $H$.
Let $D\co H\to H^*$ denote the symplectic duality map (this is the same as the Poincar\'e duality if we identify $H^*$ with $H^1(\Sigma)$). 
We have a map
 \[D\otimes id\co H\otimes\lv{k}\to H^*\otimes \lv{k}\cong \Hom(H,\lv{k}).\]
We also use $d^2$ to denote the differential \[d^2\co H_3(\gam{k})\to H_1(\gam{k},H_1(\lv{k+1}))\cong H\otimes\lv{k+1}\] of the spectral sequence of Equation~(\ref{eq:groupsequence}).

The following is a theorem of Morita~\cite[Theorem~3.1]{Morita93}.
\begin{theorem}
\label{th:MoritaTheorem}
We have:
\[\johom{k}=(D\otimes id)\circ d^2 \circ \mohom{k}\]
\end{theorem}

Let $T_{k+1}$ be the fiber of $[e]$ with respect to the map $X_{k+1}\to X_{k}$.
This space is easily seen to be a torus with fundamental group $\lv{k+1}$.
According to the theory of Section~\ref{ss:Liealgebrahomology},
we have a map $v\co C_*(T_{k+1};\R)\to C_*(\lal{k})$ inducing an isomorphism on homology.
\begin{lemma}~\label{le:commutingdifferentials}
There are homotopy equivalences 
\begin{gather*}
f\co C_*(\gam{k})\to C_*(X_{k})\\
\tag*{\text{and}} f'\co C_*(\lv{k+1})\to C_*(T_{k+1})
\end{gather*}
such that the following diagram commutes: 
\[
\xymatrix{
H_3(\gam{k})\ar[d]^{(v\circ f)_*}\ar[r]^{\hspace{-2 em}d^2} & H_1(\gam{k},H_1(\lv{k+1}))\ar[r]^{\hspace{-0.5 em}\cong} & H_1(\gam{k})\otimes H_1(\lv{k+1}) \ar[d]^{(v\circ f)_*\otimes (v\circ f')_*} \\
H_3(\lag{k})\ar[r]^{\hspace{-2 em}d^2}	      & H_1(\lag{k},H_1(\lal{k+1}))\ar[r]^{\hspace{-0.5 em}\cong} & H_1(\lag{k})\otimes H_1(\lal{k+1}), 
}
\]
where the maps $d^2$ are the differentials from the Hochschild--Serre spectral sequences for the exact sequences in Equation~(\ref{eq:groupsequence}) and Equation~(\ref{eq:algebrasequence}), respectively. 
\end{lemma}
\begin{proof}
We note that since the extensions in Equation~(\ref{eq:groupsequence}) and Equation~(\ref{eq:algebrasequence}) are central, we can make the following identifications:
\begin{align*}
H_1(\lag{k};H_1(\lal{k+1}))&=H_1(\lag{k})\otimes H_1(\lal{k+1})\\
\tag*{\text{and}}
H_1(\gam{k};H_1(\lv{k+1}))&=H_1(\gam{k})\otimes H_1(\lv{k+1})
\end{align*}
Consider the double complex:
\[C_{pq}=F_p(\gam{k})\otimes_{\gam{k}}(F_q(\gam{k+1})_{\lv{k+1}})\]
Here subscripts denote co-invariants and the boundary maps are induced from the usual ones for the bar resolutions $F_*(\gam{k})$ and $F_*(\gam{k+1})$.
This double complex gives rise to the Hochschild--Serre spectral sequence for Equation~(\ref{eq:groupsequence}).
See Brown~\cite[chapter VII.6]{Brown} for proof.

Let $U{\lag{k}}$ be the universal enveloping algebra of $\lag{k}$.
We consider a second double complex:
\[D_{pq}=(U{\lag{k}}\otimes\Lambda^p\lag{k})\otimes_{U\lag{k}}((U{\lag{k+1}}\otimes\Lambda^q\lag{k+1})_{\lal{k+1}})\]
Again, subscripts denote co-invariants.
The boundary maps here are induced from the usual ones for the Chevalley--Eilenberg complexes $U{\lag{k}}\otimes\Lambda^p\lag{k}$ and $U{\lag{k+1}}\otimes\Lambda^q\lag{k}$.
We note that the co-invariants of the Chevalley--Eilenberg complex are the standard Lie algebra chain complex with the usual boundary.
See Weibel~\cite[chapter 7.7]{Weibel} for definitions.
This double complex gives rise to the Hochschild--Serre spectral sequence for Equation~(\ref{eq:algebrasequence}).
See Weibel~\cite[chapter 7.5]{Weibel} for proof.

Now let 
$h_k\co C_*(\gam{k})\to C_*(X_k)$
be a homotopy equivalence (these exist by standard theory; see the proof to Lemma~\ref{le:reallyequivalence}). 
We also demand that $h_{k+1}|_{C_*(\lv{k+1})}$ has its image in $C_*(T_{k+1})$.
We define a map
\[f_{pq}\co C_{pq}\to D_{pq}\]
by its action on a $\Z$--basis for $C_{pq}$:
\[
\begin{split}
f_{pq}\co &[(\gamma_0,\ldots,\gamma_p)\otimes(\gamma_0',\ldots,\gamma_q')]\\
&\quad \mapsto [(1\otimes v (h_k([\gamma_0,\ldots,\gamma_p])))\otimes(1\otimes v(h_{k+1}([\gamma_0',\ldots,\gamma_q'])))]
\end{split}
\]
Here the square brackets denote the equivalence class represented by that element after taking co-invariants. 
The map $f_{**}$ is well-defined because we pass to co-invariants before plugging into $h_k$ and $h_{k+1}$.
We note that $f_{**}$ is a map of double complexes because $h_n$ and $v$ are chain maps and because the Chevalley--Eilenberg complex boundary map induces the usual Lie algebra chain boundary map.
This means that $f_{**}$ induces a map of spectral sequences and therefore:
\[d^2\circ f^2_{30}=f^2_{11}\circ d^2\co H_3(\gam{k})\to H_1(\lag{k})\otimes H_1(\lal{k+1})\]
One can easily show that
\[(v\circ h_k)_*=f^2_{30}\co H_3(\gam{k})\to H_3(\lag{k})\]
and that:
\[(v\circ h_k)_* \otimes(v\circ h_{k+1}|_{C_*(\lv{k+1})})_*=f^2_{11}\co H_1(\gam{k})\otimes H_1(\lv{k+1})\to H_1(\lag{k})\otimes H_1(\lal{k+1})\]
This proves the lemma.
\end{proof}

\begin{proof}[Proof of Main Theorem~\ref{mt:ExtJohnsonHom}]
First of all, we know from Lemma~\ref{le:d2equivariant} that the map $\exdiff\circ\emomap{k}$ is in fact a crossed homomorphism, since it is the composition of a crossed homomorphism and an equivariant map.

Let the map $\Hom(H,\lv{k+1})\into C_2(\lag{k+1})/(\lag{k+1}^{(1)}\wedge\lal{k+1})$ 
be the following long composition of the maps.
First is the symplectic duality isomorphism:
\[D^{-1}\otimes id\co \Hom(H,\lv{k+1})\isomarrow H\otimes \lv{k+1}\]
Second is this isomorphism to a product of homology groups:
\[H\otimes \lv{k+1} \isomarrow H_1(X_k)\otimes H_1(T_{k+1})\]
Then we have this injection, from Lemma~\ref{le:homologyinjection}: 
\[H_1(X_k)\otimes H_1(T_{k+1})\into H_1(X_k;\R)\otimes H_1(T_{k+1};\R)\]
This is followed by the dual Nomizu map on each part, which is an isomorphism by Corollary~\ref{co:dualofNomizuTheorem}:
\[v_*\otimes v_*\co H_1(X_k;\R)\otimes H_1(T_{k+1};\R)\isomarrow H_1(\lag{k})\otimes H_1(\lal{k+1})\]
This is then isomorphic to the following vector space:
\[H_1(\lag{k})\otimes H_1(\lal{k+1})\isomarrow (\lag{k+1}\wedge\lal{k+1})/(\lag{k+1}^{(1)}\wedge\lal{k+1})\]
This vector space then injects in our desired target:
\[(\lag{k+1}\wedge\lal{k+1})/(\lag{k+1}^{(1)}\wedge\lal{k+1})\into C_2(\lag{k+1})/(\lag{k+1}^{(1)}\wedge\lal{k+1})\]

\begin{figure}
\begin{center}
\[
\xymatrix{
\bdI{k}\ar[rr]^{\johom{k}}\ar[dr]^{\mohom{k}}\ar[ddr]_{\emomap{k}|_{\bdI{k}}} & &\Hom(H,\lv{k+1}) \ar[d]^{D^{-1}\otimes id}\\
 	& H_3(\gam{k})\ar[r]^{d^2}\ar@{^{(}->}[d] & H\otimes \lv{k+1}\ar@{^{(}->}[d] \\
	& H_3(\lag{k})\ar[r]^{\hspace{-2.5 em}d^2} & H_1(\lag{k})\otimes H_1(\lal{k+1})
}
\]
\caption{The diagram relating $d^2\circ\emomap{k}$ and $\johom{k}$.}
\label{fi:diagram}
\end{center}
\end{figure}

Now consider the diagram in Figure~\ref{fi:diagram}.
The top quadrilateral commutes by Theorem~\ref{th:MoritaTheorem}, the side triangle commutes by Main Theorem~\ref{mt:ExtMoritaHom} and the lower right square commutes by Lemma~\ref{le:commutingdifferentials}.
By Lemma~\ref{le:topologytoalgebra} we know that $\emomap{k}$ takes values in $H_3(\gam{k})$.
Since we can see that $\exdiff\circ\emomap{k}$ restricts to $d^2\circ \emomap{k}$ on $\bdI{k}$, the above commuting diagram proves the theorem.
\end{proof}

For the proof of Main Corollary~\ref{mc:sdproduct}, we need the following result of A. Heap from~\cite[Corollary~23]{Heap}.
\begin{theorem}[Heap]
The kernel of $\mohom{k}$ is $\bdI{2k-1}$.
\end{theorem}
\begin{proof}[Proof of Main Corollary~\ref{mc:sdproduct}]
Since $\bdI{k}$ acts trivially on $\gam{k}$ and the other actions under consideration are induced from this action, we see that $\bdI{k}$ is in the kernel of the actions of $\bdM$ on $\CmodB{k}$ and $C_2(\lag{k+1})/(\lag{k+1}^{(1)}\wedge\lal{k+1})$.
This means that the semi-direct products in the statement of the corollary are well-defined.
The map 
\begin{align*}
\bdM&\into(\bdM/\bdI{k})\ltimes (\CmodB{k})\\
\tag*{\text{by}} \alpha&\mapsto(\alpha\bdI{k},\emomap{k}(\alpha))
\end{align*}
is a homomorphism because the actions in the definitions of $\emomap{k}$ and the semi-direct product match.
The kernel of this map is the kernel of $\mohom{k}$: if $\alpha$ maps to the trivial element, it must be in $\bdI{k}$, on which $\emomap{k}$ is equal to $\mohom{k}$.
Since the kernel of $\mohom{k}$ is $\bdI{2k-1}$, the first map in the statement of the theorem is the induced map and is a well-defined injection.
The proof that the second map is a well-defined injection is similar.
\end{proof}

\section{Final remarks}\label{se:finalremarks}
Many questions arise concerning the best possible range for a crossed homomorphism extending a Johnson or Morita homomorphism to the mapping class group.
In the cases that have been computed explicitly, by Morita in~\cite{Morita93Invent}, and by Perron in~\cite{Perron}, there are crossed homomorphisms extending $\johomg{2}$ and $\johomg{3}$, with ranges in the finite rank abelian groups $\frac{1}{2}\Hom(H,\lv{2})$ and $\frac{1}{8}\Hom(H,\lv{3})$, respectively.
Although it is conceivable that the methods in this paper could be modified to give similar results for $\mohom{k}$ for arbitrary $k$, there are two difficulties to be overcome: 
\begin{itemize}
\item the use of integration in these methods, and the choice of $i:\Sigma\to X_k$ would have to be carefully managed to ensure well-controlled, rational coefficients, and
\item One can easily prove from the definitions that $\partial\circ \emomap{k}$ is a refinement of $\exdiff\circ\emomap{k-1}$.
This indicates that no choice of $i$ will make $\emomap{k}$ take values in $H_3(\lag{k})$ for $k>2$.
In order to restrict the range of $\emomap{k}$ from $C_3(\lag{k})/B_3(\lag{k})$ to $H_3(\lag{k})$, we will need to correct $\emomap{k}$ by a cocycle 
that projects to $\emomap{k-1}$ but restricts to the zero map on $\bdI{k}$.  It is not at all clear whether such a cocycle should exist for $k>3$.
\end{itemize}

There is an advantage to the use of integration in these methods; namely, that it leads to connections between algebraic invariants such as the Morita homomorphism and differential geometric invariants such as flux.
In particular, we have used these methods to prove a pair of theorems connecting the extended flux homomorphism to $\emomap{2}$ on a surface with an area form.
We hope to present these results in a future paper.
Of course, there may be further connections between such differential invariants and $\emomap{k}$ for higher $k$.

We note that the methods of this paper also work for the mapping class group of a closed surface $\Sigma_{g,*}$ with a basepoint.
The nilmanifolds analogous to the $X_k$ in this case are the higher Albanese manifolds of a Riemann surface, and each naturally has the structure of a $\C$--manifold.
This is due to R. Hain and S. Zucker; see Hain~\cite[Definition~1.4]{Hain87} and~\cite[Remarks~1.5]{Hain87}, or Hain and Zucker~\cite[Section~5]{Hain&Zucker}. 
There may be a version of $\emomap{k}$ over a closed surface that takes advantage of this additional structure.
In the case $k=2$, Hain gives a construction in~\cite{Hain93} that seems to be related.

Finally, one can interpret a cocycle on the mapping class group as a topological cocycle in an appropriate local coefficient system on the moduli space of surfaces.
R. Penner and S. Morita do this for Morita's extension of $\johom{2}$ in~\cite{Morita&Penner}, but their methods do not seem to give a good interpretation for the cocycles presented in this paper.
It seems natural to ask how one would develop the objects of this paper from such a perspective.

%%%%%%%%%%%%%%%%%%%%   End of main body of article
%
%                             References
%
%   BiBTeX users uncomment the following line:
%
%\bibliographystyle{gtart}
%

\noindent Department of Mathematics\\
The University of Chicago\\
5734 South University Avenue\\
Chicago, IL 60637\\
E-mail: {\tt mattday@math.uchicago.edu}

\end{document}